\documentclass[12pt]{article}
\usepackage{amssymb,amsthm,latexsym,amsmath}
\usepackage{msc}
\usepackage{graphicx}
\newtheorem{definition}{Definition}[section]

\newtheorem*{theorem A}{Theorem A}
\newtheorem{lemma}{Lemma}[section]
\newtheorem{proposition}{Proposition}[section]
\newtheorem{corollary}{Corollary}[section]
\def \AT  {{\mathcal{AT}}}
\def \DT  {{\mathcal{DT}}}
\def \ATo {\stackrel {\rightarrow}{\mathcal{AT}}}
\def \DTo {\stackrel {\rightarrow}{\mathcal{DT}}}
\def \T  {{\mathcal{T}}}
\begin{document}
%
%
\title{Arrowhead and Diamond Diameters\footnote{Submitted to {\em Elec. J. Graph Theory \& Appl.}: February 22, 2020.}}

\author{Dominique D\'{e}s\'{e}rable\\
Institut National des Sciences Appliqu\'{e}es, Rennes, France\footnote{Until 2013.}\\ email:\,\texttt{domidese@gmail.com}}

\date{April 27, 2022}

\maketitle
 %
%
\begin{abstract}
\noindent {\em Arrowhead} and {\em diamond} are two hierarchical Cayley graphs defined on the triangular grid. In their undirected version, they are isomorphic and merely define two distinct representations of the same graph. This paper gives the expression of their diameter, in the oriented and non--oriented case. It also displays the full distribution of antipodals.
\end{abstract}
%
%
%
%
%
%
\section{Introduction}
\label{section:Introduction}
%
%
\noindent Networks are represented by graphs or digraphs: a vertex in the graph stands for a node in the network while an edge (resp.\,an arc) stands for a full--duplex (resp.\,half--duplex) communication link.
The paper is concerned with a family of tori (i.e. {\it arrowhead} and {\it diamond}) which was defined on the triangular (or ``hexavalent'') grid
\cite {Deserable-1999}.
Their related interconnection networks have several important advantages. They have a bounded degree and the highest allowed degree for a two--dimensional regular grid.
They are good hosts for embedding the normal grid. As hierarchical Cayley graphs, they allow recursive constructions and divide--and--conquer schemes for information dissemination like broadcasting and gossiping
\cite {Deserable-1997,Heydemann:Marlin:Perennes-2001}.
They could accordingly provide a promising topology for Network--on--Chip ``NoC'' architectures, once the ``diagonal'' link be suitable for on--die implementation
\cite {Jayasimha:Zafar:Hoskote-2006,Balfour:Dally-2006,Salminen:Kulmala-Hamalainen-2007,Deewakar:Chatterjee-2016}.

In general, the topologies related to plane tessellations belong to the family of multi--loop networks
\cite{Bermond:Comellas:Hsu-95}
and the hexagonal case was investigated by Morillo {\em et al.}
\cite{Morillo:Comellas:Fiol-1986}
in order to exhibit graphs with minimum diameter. They proved that the maximum order of a triple loop graph with diameter $D$ is
$
  N = 3D^2 + 3D + 1.
$
The grid representation of the graph is a hexagonal torus with $D$ circular rings of length $6D$ arranged around a central node. Incidentally, this family of ``honeycombs'' $H_n$ was encountered elsewhere, arising in various projects such as FAIM--1
\cite{Davis:Robison-1985},
Mayfly
\cite{Davis-1992},
HARTS
\cite{Chen:Shin:Kandlur-1990}
and more recently with the EJ 
\cite{Huber-1994}
networks
\cite{Albader:Bose:Flahive-2012}.

The topology of our arrowhead family
\cite {Deserable-1999}
is quite different. The construction follows a recursive scheme and yields various representations of (directed) digraphs or (undirected) graphs:
(Sierpi{\'n}ski--like
\cite {Sierpinski-1916},
hexagonal) {\em arrowheads} or (lozenged, orthogonal) {\em diamonds}.
The construction of arrowheads and diamonds, not isomorphic as digraphs, is induced by two possible orientations in the hexavalent lattice.
In their undirected version, they are isomorphic and merely define two distinct representations of the same graph.
So far, the hexagonal arrowhead underlies a cellular automata network
\cite{Deserable-2002,Deserable-2011,Deserable-2012,Deserable-2014}
whereas the orthogonal diamond, named ``$T_n$'' therein, is the subject of several works on cellular multiagent systems
\cite{Ediger:Hoffmann:Deserable-2010,Ediger:Hoffmann:Deserable-2011,Hoffmann:Deserable-2014,Hoffmann:Deserable-2015}:
performances are compared with a subnetwork of $T_n$ which is just the $2^n$--ary $2$--cube
\cite{Dally:Seitz-1986};
note that another family of ``augmented'' $k$--ary $2$--cubes was investigated elsewhere
\cite{Xiang:Stewart-2011}
for any $k$ but which coincide with $T_n$ only when $k=2^n$.

Oriented and non--oriented diameters are important parameters that define the maximum distance from any vertex to another in digraphs and graphs and they provide lower bounds for routing and global communications. This paper is devoted to their study in the arrowhead family and yields their {\em exact} value as well as the full distribution of antipodals.
%
In Sect.~\ref{section:Arrowhead and Diamond},
%
we recall some general statement concerning Cayley graphs and {\em arrowhead} and {\em diamond} are redefined from
\cite {Deserable-1999}.
%
Section~\ref{section:Non--Oriented Diameter}
%
and
%
Section~\ref{section:Oriented Diameter}
%
are devoted to the analysis of the distribution of the antipodals and the computation of the diameter, respectively for the undirected and the directed version of these graphs.
%
%
\section{Arrowhead and Diamond}
\label{section:Arrowhead and Diamond}
%
%
\noindent A Cayley graph or digraph
$
  \Gamma ({\mathcal G}, {\mathcal S})
$
is constructed from a group
$
  {\mathcal G}
$
and a generating set
$
  {\mathcal S} \in   {\mathcal G} .
$
The vertex set is
$
  {\mathcal G}
$
and the edge set 
$
  {\mathcal G} \times   {\mathcal S} .
$
Cayley graphs (the same remarks hold for digraphs) are regular of degree
$
	| \, {\mathcal S} \,  | .
$
Their edge--connectivity
$
	\lambda
$ 
(the minimum number of edges whose removal disconnects the graph) satisfies
$
	\lambda = | \, {\mathcal S} \,  | .
$
They are vertex--transitive (or vertex--symmetric) in the sense that for any pair 
$
	( u , u')
$ 
of vertices there exists an automorphism of
$
	\Gamma
$ 
that maps $u$ into $u'$.
%
%
\begin{definition}
\label{definition:  arrowhead and diamond}
%
%
Given the group
$
  G_n =  \mathbb{Z}_{2^n} \times  \mathbb{Z}_{2^n}  
$
with 
$
	n \in \mathbb{N}:
$
\begin{itemize}
\item We define the directed  {\em arrowhead}
$
  \ATo_n = \Gamma  \, (G_n, S^+)
$
as the digraph of $G_n$ with the generating set
$
	S^+ = (s_1, s_2, s_3 ) = ( (-1,-1), (1,0), (0,1) ) .
$
\item Let 
$
	S^- = (-s_1, -s_2, -s_3 )  
$
be the set of inverses of $S^+$ and let the generating set 
$
	S =  S^+ \cup S^-
$
now closed under inverses. We define the undirected {\em arrowhead} as 
$
  \AT_n = \Gamma  \, (G_n, S) .
$
%
%
\item We define the directed  {\em diamond}
$
  \DTo_n = \Gamma  \, (G_n, T^+)
$
as the digraph of $G_n$ with the generating set
$
	T^+ = (t_1, t_2, t_3 ) = (- s_1, s_2, s_3 ) .
$
\item Let 
$
	T^- = (-t_1, -t_2, -t_3 )  
$
be the set of inverses of $T^+$ and let the generating set 
$
	T =  T^+ \cup T^-
$
now closed under inverses. We define the undirected {\em diamond} as 
$
  \DT_n = \Gamma  \, (G_n, T) .
$
\item  $\DT_n$ and $\AT_n$  are isomorphic since $T = S$.
$\hfill \Box$
\end{itemize}
\end{definition}
%
%
 \noindent The order of these graphs is given by
$
  N = | \, G_n \, | = 4^n
$
and the number of arcs (or edges) by
$
	 \frac{1}{2} | \, S \, | |\,  G_n \, |  =  3  N .
$
In the sequel, the detection of antipodals and the computation of diameters will follow an inductive scheme. It is therefore convenient to exploit the recursive structure of the graphs as follows. \\
%
%
\begin{definition}
\label{definition:  Subgroups of G_n}
Let 
$
	0  \le k \le n
$
and let again
\begin{equation}
     G_{n} = \{  ( \, x,  y \, ) = x s_2 + y s_3 \, | \, x, y  \in \mathbb{Z}_{2^{n}} \} 
\end{equation}
be the vertex set expressed from 
the previous definition
and let
\begin{equation}
     G_{n,k} =   \mathbf{2}^k \mathbf{\cdot} \ G_{n-k} =   \{  ( \, 2^k x, 2^k y \, ) \, | \, x, y  \in \mathbb{Z}_{2^{n-k}} \} 
\label{equation:  Definition of of G_(n,k)}
\end{equation}
be the subgroup of 
$
	G_n
$ 
generated by
$
	\mathbf{2}^k S^+ =  \{ \, 2^k s \, | \, s \in S^+  \,  \}  .
$

We define  
$
  \ATo_{n,k}   \  = \Gamma  \, ( G_{n,k} \,, \mathbf{2}^k S^+)
$
isomorphic to $\ATo_{n-k}$ as the digraph of $G_{n,k}$ with the generating set
$
	 \mathbf{2}^k S^+ .
$
$\hfill \Box$
\end{definition}
%
%
\noindent As a matter of fact, this statement gives an embedding scheme of $\ATo_{n-k}$ into $\ATo_n$ with dilation $2^k$. 
A definition of digraph $\DTo_{n,k}$ and of their undirected counterpart will follow from a similar statement.

In the following, the undirected --isomorphic-- graphs $\AT_n$ and  $\DT_n$ will be denoted as $\T_n$.
%
%
\begin{figure}
\centering
\includegraphics[width=11cm]{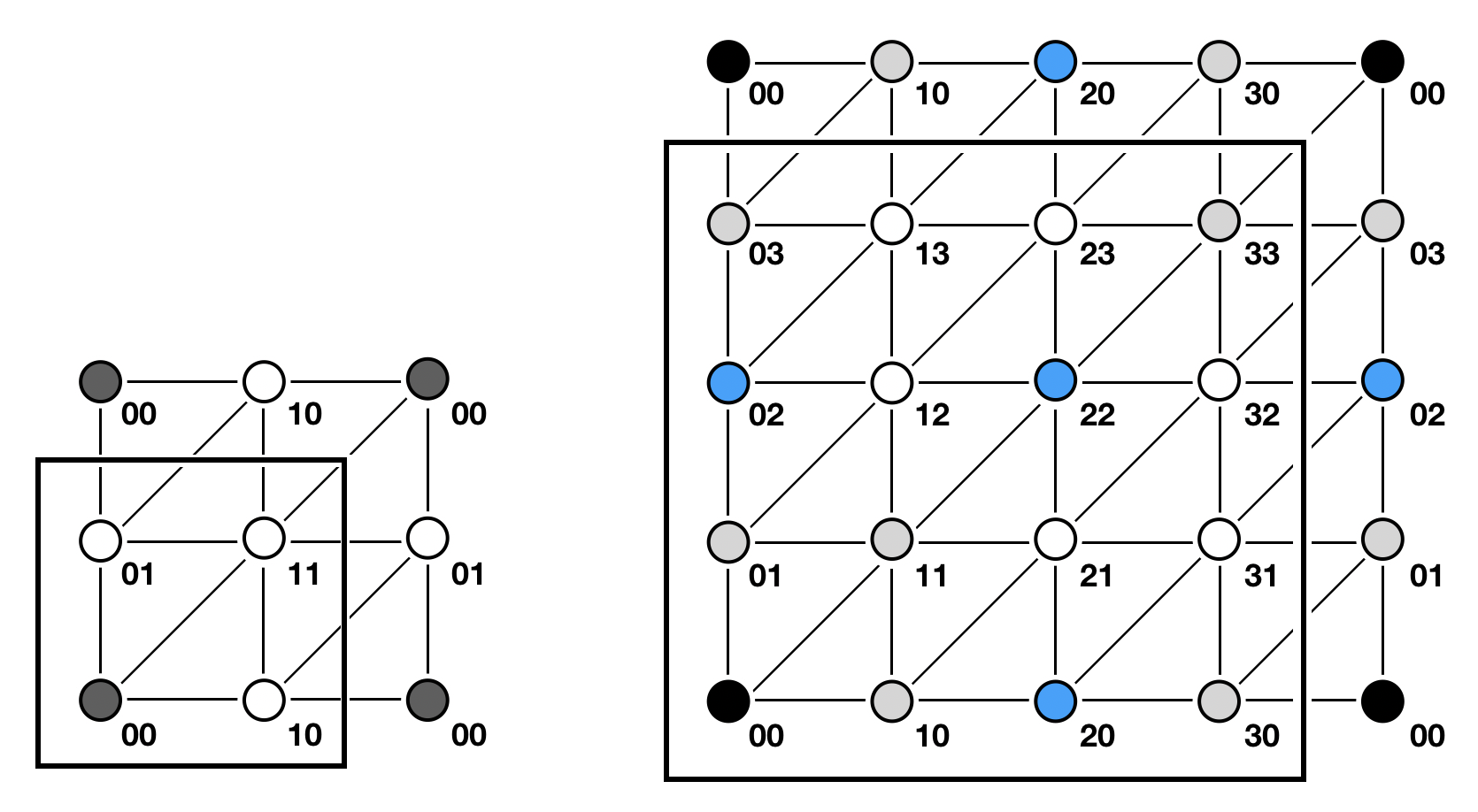}
\caption{Representation of $\T_1$ and $\T_2$ -- some vertices are replicated for convenience.}
\label{fig:Representation of T_1 and T_2}
\end{figure}
%
%
The representation of $\T_1$ and $\T_2$  is displayed in 
Fig.\,\ref{fig:Representation of T_1 and T_2}.

$\T_0$ is not displayed, reduced to the single vertex of $G_0 =   \{  ( 0,0)  \}  $ as a 6--valent vertex.
It should be observed in
(\ref{equation:  Definition of of G_(n,k)})
that
$
	 G_{n,0} =  G_{n}
$
and
$
	 G_{n,n} =  G_{0} .
$

In $\T_1$ we observe that
$
	G_{1,0} =  G_1 =   \{  ( 0,0),( 1,1),( 1,0),( 0,1) \}  
$ 
with 4 vertices and that 
$
	G_{1,1} =  \mathbf{2} \mathbf{\cdot} \ G_0  =   \{  ( 0,0)  \} .
$

In $\T_2$ we get
$
     G_{2} = \{  ( \, x,  y \, ) = x s_2 + y s_3 \, | \, x, y  \in \mathbb{Z}_4 \} 
$
with 16 vertices and see that 
$
     G_{2,1} =   \mathbf{2} \mathbf{\cdot} \ G_{1} =   \{  ( 0,0),( 2,2),( 2,0),( 0,2) \} 
$
isomorphic to $G_1$ and that
$
     G_{2,2} =   \mathbf{4} \mathbf{\cdot} \ G_{0} =    \{  ( 0,0)  \} .
$
%
%
\section{Non--Oriented Diameter}
\label{section:Non--Oriented Diameter}
%
%
\noindent The diameter of a graph is the maximum distance between any pair of vertices. We call ``oriented'' (resp.\,``non--oriented'') the diameter of a directed (resp.\,undirected) graph. Two vertices are said to be {\em antipodal} if their distance is the diameter. Without loss of generality, we can compute the diameter as the length of the shortest path from the ``origin'' $(0,0)$ to an antipodal since the graphs are vertex--transitive. Hereafter, the term ``antipodal'' will be viewed from {\em that} vertex.
%
%
\subsection{Diameter in $\T_n$}
\label{subsection:Diameter in Tn}
\begin{figure}
\centering
\includegraphics[width=12cm]{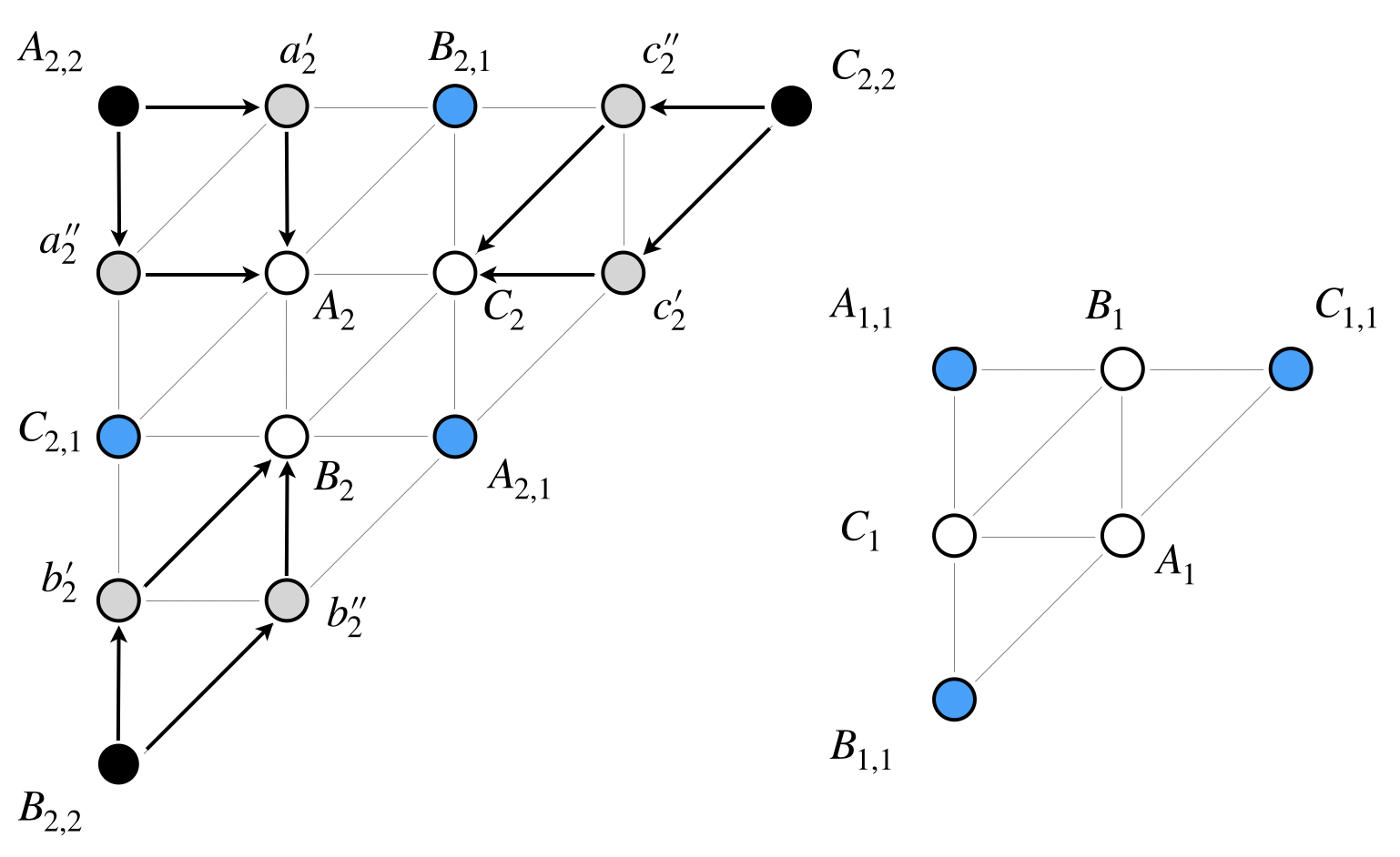}
\caption{Triangular diagram of $\T_2$ and  $\T_1$. Induction from $\T_1$ to $\T_2$. Resulting shortest paths to antipodals in $\T_2$.}
\label{fig:Triangular diagram of T_2 and T_1}
\end{figure}

%
%
\noindent In
Fig.\,\ref{fig:Representation of T_1 and T_2},
it appears  as trivial that vertices in subset 
$
	 \{ ( 1,1),( 1,0),( 0,1) \}  
$ 
in  $\T_1$ are antipodal and at distance 1.
We claim that in $\T_2$ the subset 
$
	 \{ ( 1,2),( 1,3),( 2,3) \}  
$ 
as well as its symmetric part
$
	 \{ ( 2,1),( 3,1),( 3,2) \}  
$ 
are antipodal at distance 2. To sketch the proof, let us fix 
$
	x \le y
$
and then extend to the triangular diagrams in 
Fig.\,\ref{fig:Triangular diagram of T_2 and T_1}.
Beforehand we give the following definitions.
%
%
\begin{definition}
\label{definition:  Induction on antipodal subsets}
Herein, for $n>0$, the ordered subset
$
    \Omega_{n} 
$
will always denote an antipodal 3--cycle in $\T_n$ and
$
    \Omega_{n,1} 
$
will denote the image of 
$
    \Omega_{n-1}
$
whereas
$
    \Omega_{n,2} 
$
will denote the image of 
$
    \Omega_{n-1,1}
$
under the mapping induced by
{\em (\ref{fig:Triangular diagram of T_2 and T_1})}.
\begin{itemize}
\item 
$
	\Omega_0  =   \{  ( 0,0)  \} 
$
\item  
$
	\Omega_1  = (  (1,1),( 1,0),( 0,1)  ) 
$
\item  
$
	\Omega_{n,1} =   \mathbf{2}\mathbf{\cdot}  \Omega_{n-1} =   \{  ( \, 2x, 2y \, ) \, | \, x, y  \in  \Omega_{n-1}  \}  \hspace{30mm}  (n>0)
$
\item  
$
	\Omega_{n,2} =   \mathbf{2}\mathbf{\cdot}  \Omega_{n-1,1} =   \{  ( \, 2x, 2y \, ) \, | \, x, y  \in  \Omega_{n-1,1}  \}  =  
	\mathbf{4}\mathbf{\cdot}  \Omega_{n-2} \hspace{8mm}  (n>1) 
$
$\hfill \Box$
\end{itemize}
\end{definition}
%
%
%
%
%
\noindent In the sequel, the term {\em image} will always mean the image under the isomorphic mapping in 
%
%
Def.\,\ref{definition:  Induction on antipodal subsets}.
%
%
Note that any 
$
           \Omega_{n,k}
$
is a subset of $G_{n,k}$. In the diagrams of 
Fig.\,\ref{fig:Triangular diagram of T_2 and T_1}
and further, the vertices of those different subsets will be labeled as follows:
\begin{itemize}
\item  
$
	\Omega_n  = ( A_n, B_n, C_n ) 
$
\item  
$
	\Omega_{n,1} =  ( A_{n,1}, B_{n,1}, C_{n,1} ) 
$
\item  
$
	\Omega_{n,2} =  ( A_{n,2}, B_{n,2}, C_{n,2} ).
$ 
\end{itemize}
%
%
\begin{figure}
\centering
\includegraphics[width=12cm]{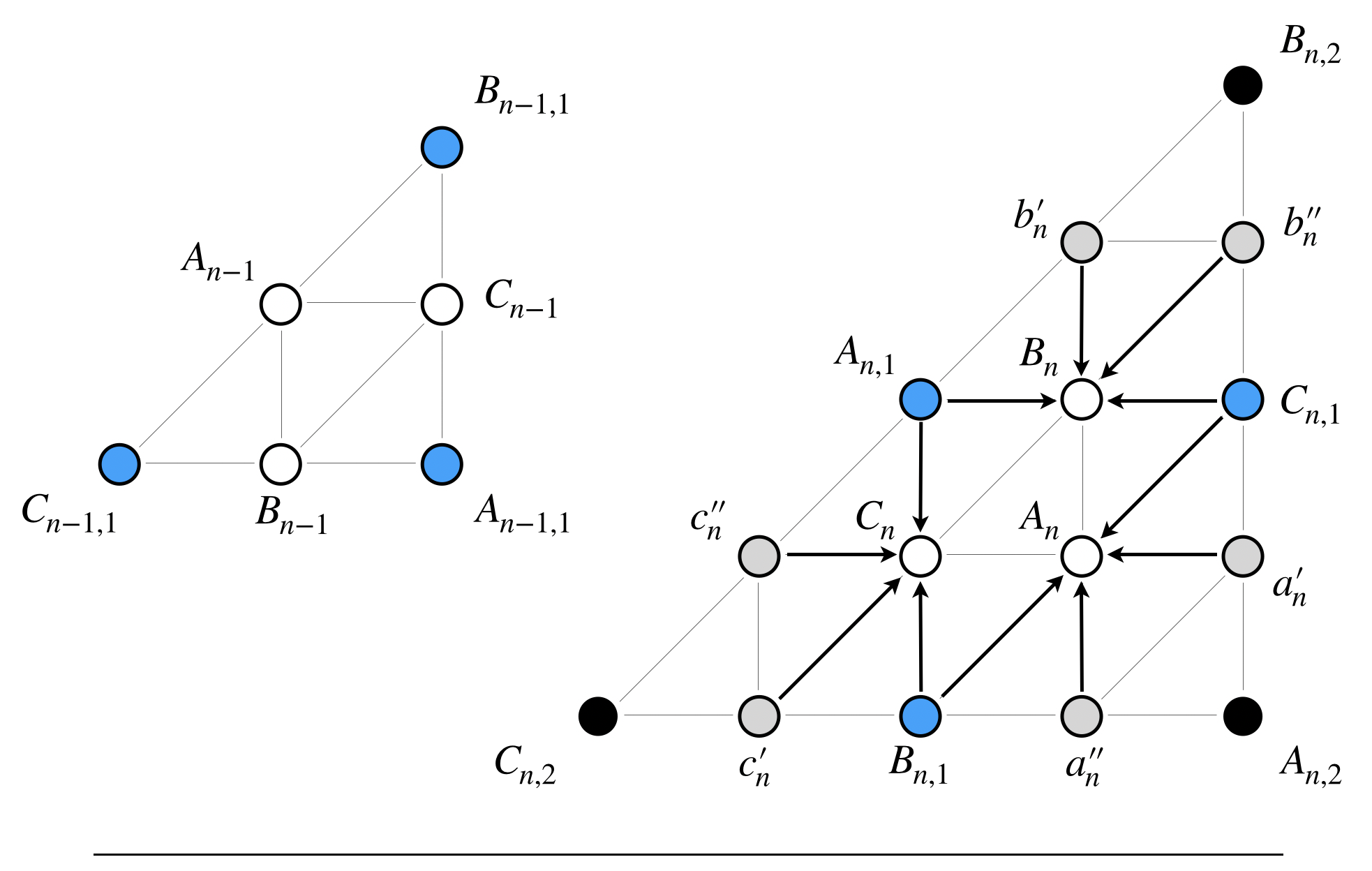}
\includegraphics[width=12cm]{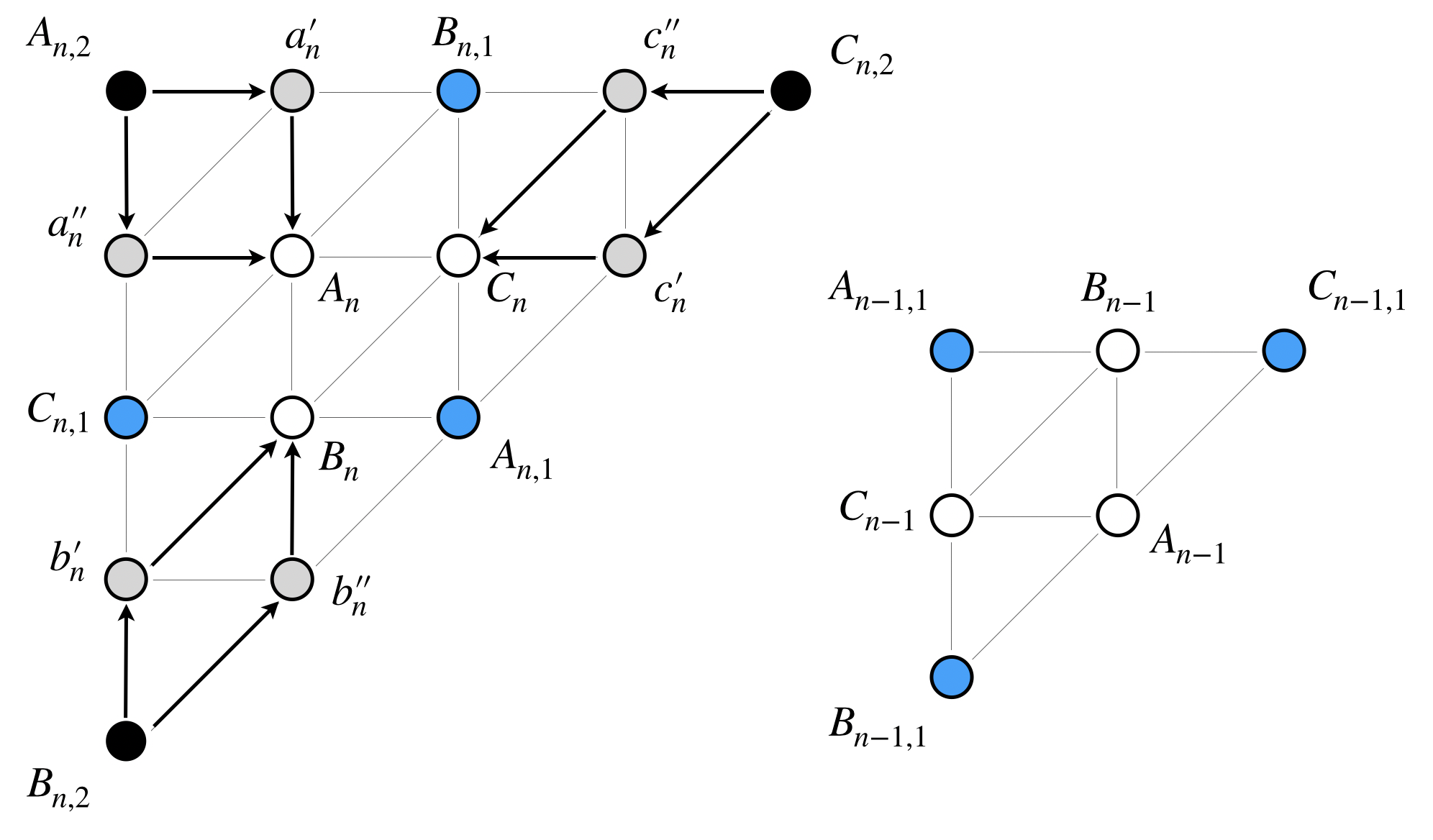}
\caption{Triangular diagrams  in $\T_n$ -- for $n$ odd $(\uparrow)$ -- for $n$ even $(\downarrow)$.}
\label{fig:Triangular diagram for n odd and n even}
\end{figure}
%
%
%
\begin{lemma}
\label{lemma:  Induction on Non--Oriented Diameter}
Let $D_n$ be the non--oriented diameter of $\T_n$.  Then $D_0 = 0$ and 
$
	D_n = 2 D_{n-1}  + 1
$
or
$
	D_n = 2 D_{n-1} 
$
depending upon the odd--even parity of n.
\end{lemma}
%
%
\textit{Proof}:
We prove by induction and refer first to
Fig.\,\ref{fig:Triangular diagram of T_2 and T_1}.
For the base case $n=1$  the vertices of 
$
	\Omega_{1,1}  
$
coincide as origin at distance 0 and vertex set 
$
	\Omega_1   
$
is trivially antipodal at distance $D_1 = 1$.
For $n=2$, vertex set
$
	\Omega_{2,2}  
$
as image of 
$
	\Omega_{1,1}  
$
coincide as origin at distance 0.
Now, vertices 
$
	( a'_2, a''_2),( b'_2, b''_2),( c'_2, c''_2)
$
reachable from
$
	( A_{2,2}, B_{2,2}, C_{2,2} ) 
$
are at distance 1 and then the triple 
$
	 ( A_2, B_2, C_2 ) 
$
is at distance 2. Besides, vertices of 
$
	\Omega_{2,1}  
$
as image of
$
	\Omega_1   
$
are by induction at distance 
$
	2 D_1 = 2.
$
Therefore the triple 
$
	\Omega_2
$
is antipodal, as well as 
$
	\Omega_{2,1}  
$
and at distance the diameter $ D_2 = 2$.
%
%
\begin{itemize}
%
%
\item  $n$  {\em odd}:                         
%
%
Let
$
	\Omega_{n,1}  
$
and
$
	\Omega_{n,2}  
$
be the respective images of
$
	\Omega_{n-1}  
$
and
$
	\Omega_{n-1,1}  
$
in upper
Fig.\,\ref{fig:Triangular diagram for n odd and n even}.
From induction, vertices of those both triples are assumed at distance 
$
	2 D_{n-1} .
$
Now the six vertices 
$
	( a'_n, a''_n),( b'_n, b''_n),( c'_n, c''_n)
$
are also at distance
$
	2 D_{n-1}:
$
$a'_n$ and $A_{n,2}$  are reachable from a common antecedent
\[
	(  \, x_{a'_n}  , y_{a'_n} \, ) + t_2 =  (  \, x_{A_{n,2}}  , y_{A_{n,2}}  \, ) + t_1
\]
and idem for $a'_n$ and $C_{n,1}$ 
\[
	(  \, x_{a'_n}  , y_{a'_n} \, ) + t_1 =  (  \, x_{C_{n,1}}  , y_{C_{n,1}}  \, ) + t_2
\]
and so forth.Therefore those six vertices are at distance
$
	2 D_{n-1} .
$
Vertex $A_n$ can then be reached through
$
	( a'_n, a''_n), 
$
$B_n$ through
$
	( b'_n, b''_n), 
$
$C_n$ through
$
	 ( c'_n, c''_n).
$
The triple
$
	( A_n, B_n, C_n ) 
$
can also be reached from
\[
	(B_{n,1}, C_{n,1}), (C_{n,1}, A_{n,1}), (A_{n,1}, B_{n,1}).
\]
It is therefore antipodal, at distance the diameter 
$
	D_n = 2 D_{n-1}  + 1 .
$
%
%
%
%
\item  $n$  {\em even}:                         
%
Let
$
	\Omega_{n,2}  
$
be the image of
$
	\Omega_{n-1,1}  
$
in lower
Fig.\,\ref{fig:Triangular diagram for n odd and n even}
and, from induction, with vertices assumed at distance 
$
	2 (D_{n-1}  - 1).
$
The three pairs 
$
	( a'_n, a''_n),( b'_n, b''_n),( c'_n, c''_n)
$
reachable from
$
	( A_{n,2}, B_{n,2}, C_{n,2} ) 
$
are at distance
$
	2 (D_{n-1}  - 1) + 1 = 2 D_{n-1}  - 1 .
$
Vertex $A_n$ can be reached through
$
	( a'_n, a''_n), 
$
$B_n$ through
$
	( b'_n, b''_n), 
$
$C_n$ through
$
	 ( c'_n, c''_n)
$
and then the triple
$
	( A_n, B_n, C_n ) 
$
is at distance 
$
	( 2 D_{n-1}  - 1) + 1 =  2 D_{n-1} .
$
Besides, 
$
	\Omega_{n,1}  
$
as image of
$
	\Omega_{n-1}  
$
has also its vertices at distance
$
	2 D_{n-1} .
$
Therefore the triple 
$
	( A_n, B_n, C_n ) 
$
is antipodal, as well as 
$
	( A_{n,1}, B_{n,1}, C_{n,1} ) 
$
and at distance the diameter
$
	D_n = 2 D_{n-1} .
$
$\hfill \Box$
\end{itemize}
%
%
%
%
\begin{lemma}
\label{lemma:  Relations on diameter}
Two useful relations are exhibited.
\begin{equation}
	  \forall n>0 :   D_{n-1} + D_n  = 2^n - 1
 \label{equation:  Relation 1 on diameter}
\end{equation}
\begin{equation}
	\forall n>1 :    D_n - D_{n-2}  =  2^{n-1}  \  \    
 \label{equation:  Relation 2 on diameter}
\end{equation}
\end{lemma}
%
%
\textit{Proof}:
\label{proof:  Relations on diameter}
Let us rewrite
$
	D_n = 2 D_{n-1} + \varepsilon_n  
$
from 
Lemma \,\ref{lemma:  Induction on Non--Oriented Diameter}
where
$
	\varepsilon_n	\equiv	n \pmod{2}.
$
Given
$
	u_n = D_{n-1} + D_n 
$
we note that $u_1 = 1$ and show that 
$
	u_{n+1}  = 2 u_n + 1.
$
Now 
$
	u_{n+1}  =  D_n + D_{n+1}  =  ( 2 D_{n-1} + \varepsilon_{n}  ) +  ( 2 D_{n} + \varepsilon_{n+1}  ) 
$
but
$
	 \varepsilon_{n}  +  \varepsilon_{n+1} = 1
$
for any $n$ whence 
$
	u_{n+1}  = 2 u_n + 1 
$
and 
(\ref{equation:  Relation 1 on diameter})
results from an obvious induction.
 
For
(\ref{equation:  Relation 2 on diameter})
we note from above that
$
	D_n - D_{n-2}  =  u_{n}  - u_{n-1}  = 2^{n-1} .
$
$\hfill \Box$
%
%
%
\begin{proposition}
\label{proposition:  Tn Diameter}
$\T_n$ has the diameter
\[
	D_n  = \frac{2 \sqrt{N}-1}{3}     \    \    \   \mbox{or}   \    \    \  D_n  = \frac{2 ( \sqrt{N}-1 ) }{3}  
\]
depending on the odd--even parity of $n$.
\end{proposition}
%
%
\noindent \textit{Proof}:
Follows from
Lemma \,\ref{lemma:  Induction on Non--Oriented Diameter}.
For $n$ odd,
$
	D_n = 2 D_{n-1} + 1
$
and by
(\ref{equation:  Relation 1 on diameter})
$
	D_{n-1}  = ( 2^n - 1 ) - D_n  
$
 whence
 $
 	3 D_n =  2 \cdot 2^{n} - 1 .
 $
For $n$ even,
$
	D_n = 2 D_{n-1} 
$
 whence
 $
 	3 D_n =  2 \cdot( 2^{n} - 1) .
 $
$\hfill \Box$  \\ 
%
%

\noindent The ten first values of $\T_n$ diameter are displayed below.
\vspace{5mm}

\begin{tabular}{|c||c|c|c|c|c|c|c|c|c|c|}
 \hline
     $n$        & 0  & 1 & 2  &  3   & 4     &     5     &    6        &    7       &     8        &      9       \\ \hline
    $N$        & 1 & 4 & 16  & 64 & 256  & 1024   & 4096    & 16384  & 65536   & 262144    \\ \hline
     $D_n$   & 0 & 1 &  2   &  5  &  10   &   21    &   42       &    85     &    170      &    341     \\ \hline
\end{tabular}
%
%
%
\subsection{Antipodals in $\T_n$}
\label{subsection: Antipodals in Tn}
\begin{figure}
\centering
\includegraphics[width=11cm]{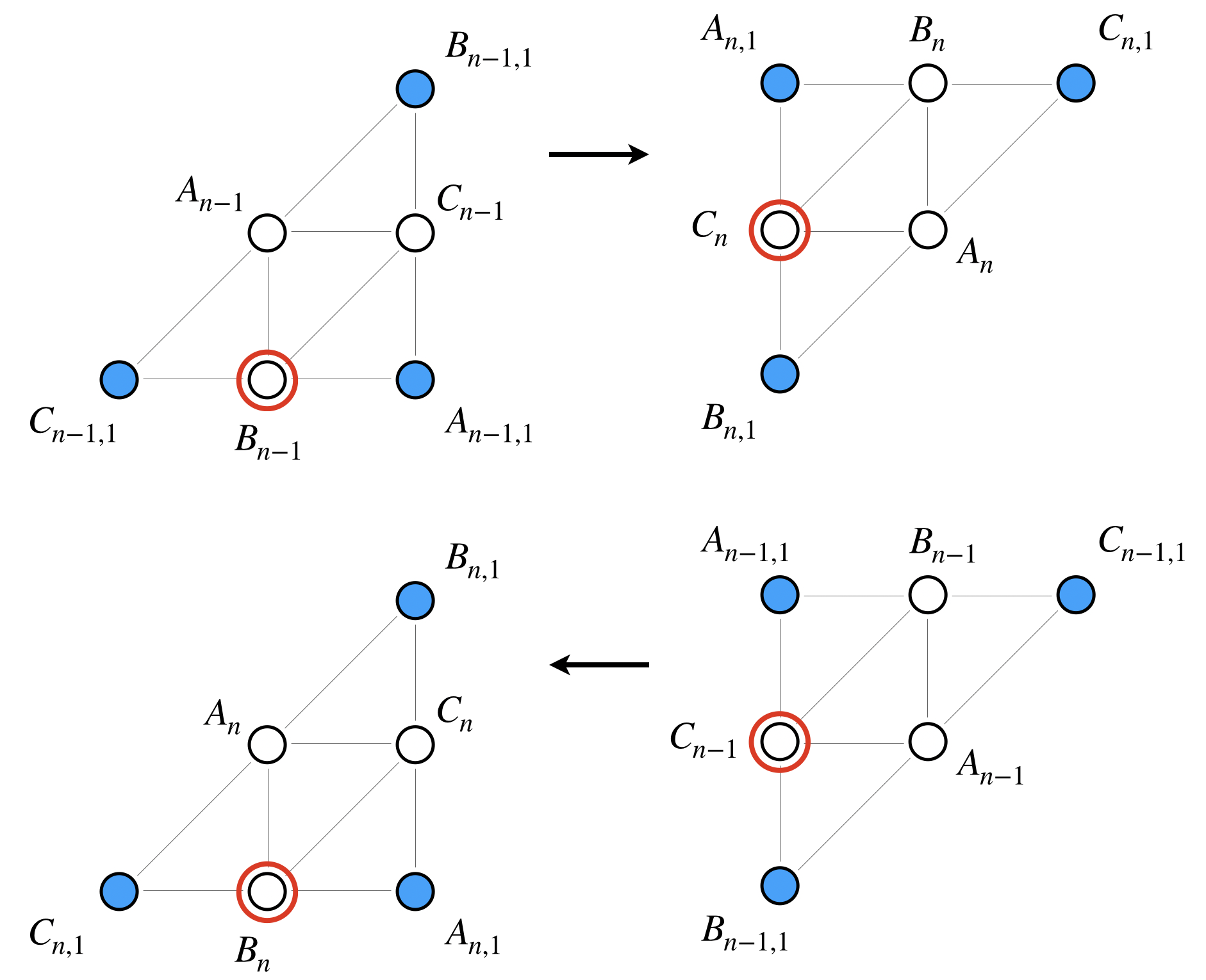}
\caption{Antipodal coordinates in  $\T_n$  -- for $n$ odd $(\uparrow)$ -- for $n$ even $(\downarrow)$.}
\label{fig:Antipodal coordinates}
\end{figure}

%
%
%
\noindent Some relevant properties of antipodal coordinates are highlighted hereafter and antipodals are enumerated.
%
%
%
\begin{proposition}
\label{proposition:  antipodal coordinates}
The antipodal coordinates in 
$
	{\Omega}_{n} 
$
 satisfy:
%
%
\begin{itemize}
%
%
\item  $n$  {\em odd}:                         
$
	(  \, x_{C_{n}}  , y_{C_{n}}  \, )  = (  \, D_{n-1}  , D_n  \, )   
$
%
%
\item  $n$  {\em even}:                         
$
	(  \, x_{B_{n}}  , y_{B_{n}}  \, )  = (  \, D_{n-1}  , D_n  \, ) .
$
%
%
\end{itemize}
%
%
\end{proposition}
%
%
\noindent \textit{Proof}:
This is true for $n = 1$ where
$
	(  \, x_{C_{1}}  , y_{C_{1}}  \, )  = ( 0,1) = (  \, D_{0}  , D_1  \, )
$
and for $n = 2$ where
$
	(  \, x_{B_{2}}  , y_{B_{2}}  \, )  = ( 1,2) = (  \, D_{1}  , D_2  \, )
$
referring back to
Fig.\,\ref{fig:Representation of T_1 and T_2}
and
Fig.\,\ref{fig:Triangular diagram of T_2 and T_1}.
%
%
%
\begin{itemize}
%
%
\item  $n$  {\em odd}     
(upper Fig.\,\ref{fig:Antipodal coordinates}):
assume
$
	(  \, x_{B_{n-1}}  , y_{B_{n-1}}  \, )  = (  \, D_{n-2}  , D_{n-1}   \, ) .
$

Then 
$
	(  \, x_{B_{n,1}}  , y_{B_{n,1}}  \, )  = (  \, 2 D_{n-2}  , 2 D_{n-1}   \, ) = (  \, D_{n-1}  , D_n  \, -1 )   
$
since $n-1$ is even when $n$ is odd, whence
$
	(  \, x_{C_{n}}  , y_{C_{n}}  \, )  =  (  \, x_{B_{n,1}}  , y_{B_{n,1}}  \, )  + (0,1).  
$
%
%
\item  $n$  {\em even}                       
(lower Fig.\,\ref{fig:Antipodal coordinates}):
assume
$
	(  \, x_{C_{n-1}}  , y_{C_{n-1}}  \, )  = (  \, D_{n-2}  , D_{n-1}   \, ) .
$

Then 
$
	(  \, x_{C_{n,1}}  , y_{C_{n,1}}  \, )  = (  \, 2 D_{n-2}  , 2 D_{n-1}   \, ) = (  \, D_{n-1} -1 , D_n  )   
$
since $n-1$ is odd when $n$ is even, whence
$
	(  \, x_{B_{n}}  , y_{B_{n}}  \, )  =  (  \, x_{C_{n,1}}  , y_{C_{n,1}}  \, )  + (1,0).  
$
$\hfill \Box$   
%
%
\end{itemize}
%
%
Coordinates of other antipodals are simply deduced.
%
%
%
%
%

\begin{figure}
\centering
\includegraphics[width=11cm]{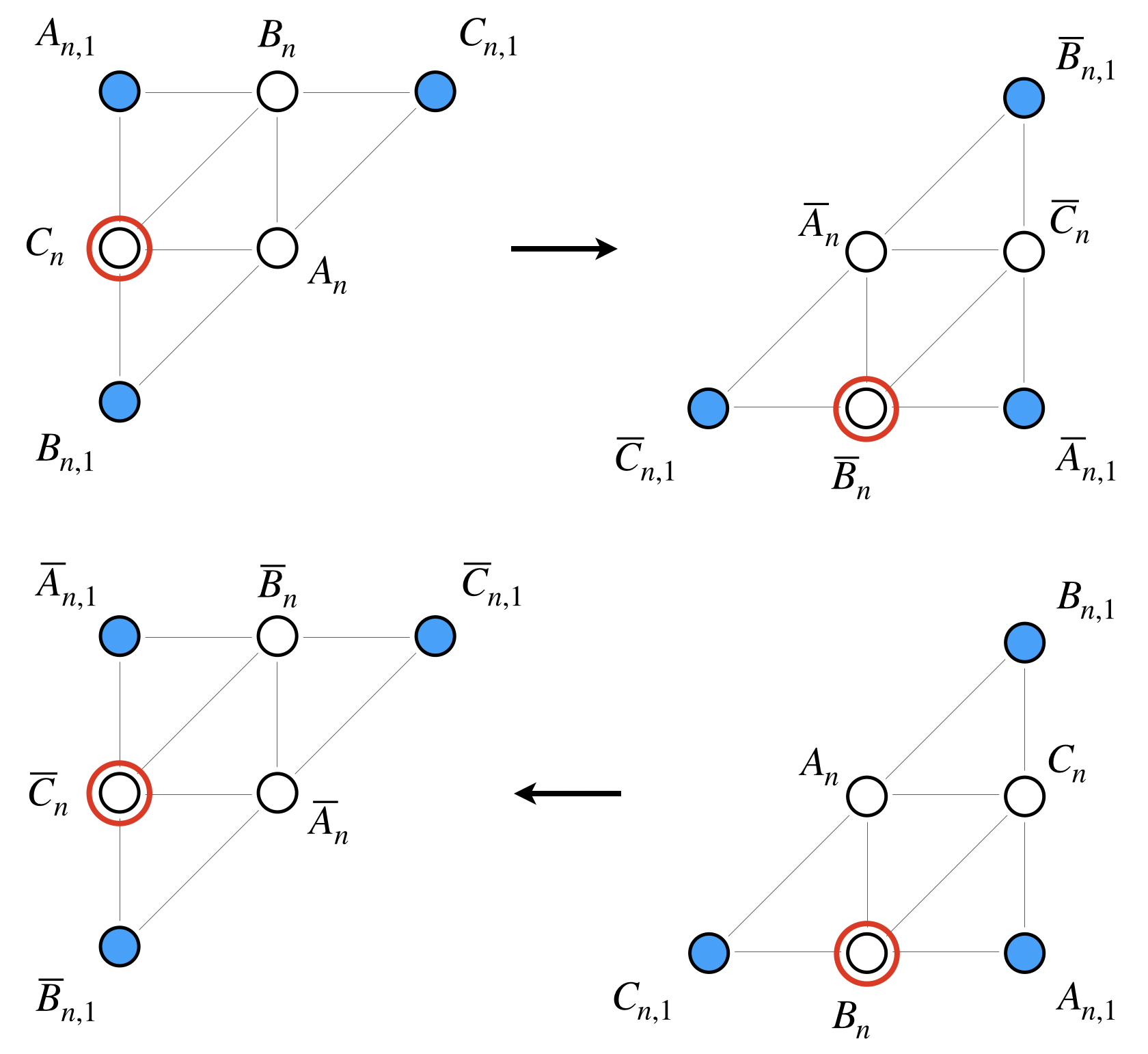}
\caption{Antipodal inverses in  $\T_n$   -- for $n$ odd $(\uparrow)$ -- for $n$ even $(\downarrow)$.}
\label{fig:Antipodal inverses}
\end{figure}

Referring back to 
Fig.\,\ref{fig:Representation of T_1 and T_2}
we now focus on the subset
$
	 \{ ( 2,1),( 3,1),( 3,2) \} , 
$ 
antipodal by  {\em symmetry}. By fixing
$
	x \ge y 
$
and as from 
%
Def.\,\ref{definition:  Induction on antipodal subsets}
%
we define the symmetric subsets
\begin{itemize}
\item  
$
	 \overline{\Omega}_{n} = ( \overline{A}_{n}, \overline{B}_{n}, \overline{C}_{n} ) 
$
\item  
$
	 \overline{\Omega}_{n,1} = ( \overline{A}_{n,1}, \overline{B}_{n,1}, \overline{C}_{n,1} ) 
$
\item  
$
	 \overline{\Omega}_{n,2} = ( \overline{A}_{n,2}, \overline{B}_{n,2}, \overline{C}_{n,2} ) 
$
\end{itemize}
where
$
	(  \, x_{ \overline{A}_{n}}  , y_{ \overline{A}_{n}}  \, )  = (  \,- x_{A_{n}}  ,- y_{A_{n}}  \, )  
$
and where any vertex in those subsets is defined in that way.
%
%
\begin{proposition}
\label{proposition:  antipodal inverses}
The antipodal coordinates in 
$
	 \overline{\Omega}_{n} 
$
 satisfy:
%
%
\begin{itemize}
%
%
\item  $n$  {\em odd}:                         
$
	(  \, x_{ \overline{B}_{n}}  , y_{ \overline{B}_{n}}  \, )  = (  \,  D_n , D_{n-1}  \, )   
$
%
%
\item  $n$  {\em even}:                         
$
	(  \, x_{ \overline{C}_{n}}  , y_{ \overline{C}_{n}}  \, )  = (  \,  D_n , D_{n-1}  \, ).
$
%
%
\end{itemize}
%
%
\end{proposition}
%
%
\noindent \textit{Proof}:
From Proposition 
\ref{proposition:  antipodal coordinates} 
and
Fig.\,\ref{fig:Antipodal inverses}.
%
%
\begin{itemize}
%
%
\item  $n$  {\em odd}:   
$
	(  \, x_{B_{n}}  , y_{B_{n}}  \, )  =  (  \, x_{C_{n}}  , y_{C_{n}}  \, )  + (1,1) 
$
whence 
\[
	(  \, x_{ \overline{B}_{n}}  , y_{ \overline{B}_{n}}  \, )  = (  \, - D_{n-1} - 1 , - D_n - 1 \, )
\]
and 
(\ref{equation:  Relation 1 on diameter})
yields
$
	( \, x_{ \overline{B}_{n}}  , y_{ \overline{B}_{n}} \,)  = ( \,D_n , D_{n-1} \,)
$
by a simple reduction in
$
\mathbb{Z}_{2^n}.
$
%
%
\item  $n$  {\em even}:                       
$
	(  \, x_{C_{n}}  , y_{C_{n}}  \, )  =  (  \, x_{B_{n}}  , y_{B_{n}}  \, )  + (1,1) 
$
whence 
\[
	(  \, x_{ \overline{C}_{n}}  , y_{ \overline{C}_{n}}  \, )  = (  \, - D_{n-1} - 1 , - D_n - 1 \, )
\]
whence again the result.
$\hfill \Box$   
%
%
\end{itemize}
%
%
%
%
%
%
%
\begin{corollary}
\label{corollary:  antipodal inverses}
The antipodal inverses satisfy for any $n$:
%
%
\begin{itemize}
%
%
\item                  
$
	(  \, x_{ \overline{A}_{n}}  , y_{ \overline{A}_{n}}  \, )  = (  \, y_{A_{n}}  , x_{A_{n}}  \, )  
$
%
%
\item                 
$
	(  \, x_{ \overline{B}_{n}}  , y_{ \overline{B}_{n}}  \, )  = (  \, y_{C_{n}}  , x_{C_{n}}  \, )  
$
%
%
\item                        
$
	(  \, x_{ \overline{C}_{n}}  , y_{ \overline{C}_{n}}  \, )  = (  \, y_{B_{n}}  , x_{B_{n}}  \, )  
$
$\hfill \Box$
%
%
%
\end{itemize}
%
%
\end{corollary}
%
%
%
%
Coordinates of other antipodals subsets are simply deduced.
%
%
\begin{corollary}
\label{Corollary: Number of Antipodals in T_n}
Let
$
    {\mathcal{N_{A}}}_{,n}
$
be the number of antipodals in $\T_n$, then
%
%
\begin {itemize}
\item
$
          {\mathcal{N_{A}}}_{,1}$ = $ \mid \Omega_{1} \cup \overline{\Omega}_{1} \mid \ = \ \mid \Omega_{1} \mid \ = 3
$
\item
%
    ${\mathcal{N_{A}}}_{,2}$ = $ \mid \Omega_{2} \cup \overline{\Omega}_{2} \mid
                         +   \mid \Omega_{2,1} \cup  \overline{\Omega}_{2,1}  \mid \
                         = \ \mid \Omega_{2} \cup  \overline{\Omega}_{2}  \mid + \mid \Omega_{2,1} \mid \ = 9$
\item
$
    {\mathcal{N_{A}}}_{,n} =
           \left\{ \begin{array}{ll}
                        \mid \Omega_{n} \cup \overline{\Omega}_{n} \mid \ = \ 6 \\
                        \mid \Omega_{n} \cup  \overline{\Omega}_{n} \mid
                    +   \mid \Omega_{n,1} \cup  \overline{\Omega}_{n,1}  \mid \ = 12
                   \end{array}
           \right.
$
\end {itemize}
%
%
depending on the odd--even parity of $n > 2$.
$\hfill \Box$
\end{corollary}
%
%
\noindent We observe that 
$
	\Omega_{1} 
$
and
$
	\overline{\Omega}_{1}
$
coincide and that
$
	\Omega_{2,1}
$
and
$
	 \overline{\Omega}_{2,1} 
$
coincide.
%
%
\section{Oriented Diameter}
\label{section:Oriented Diameter}
%
%
%
\begin{figure}
\centering
\includegraphics[width=12cm]{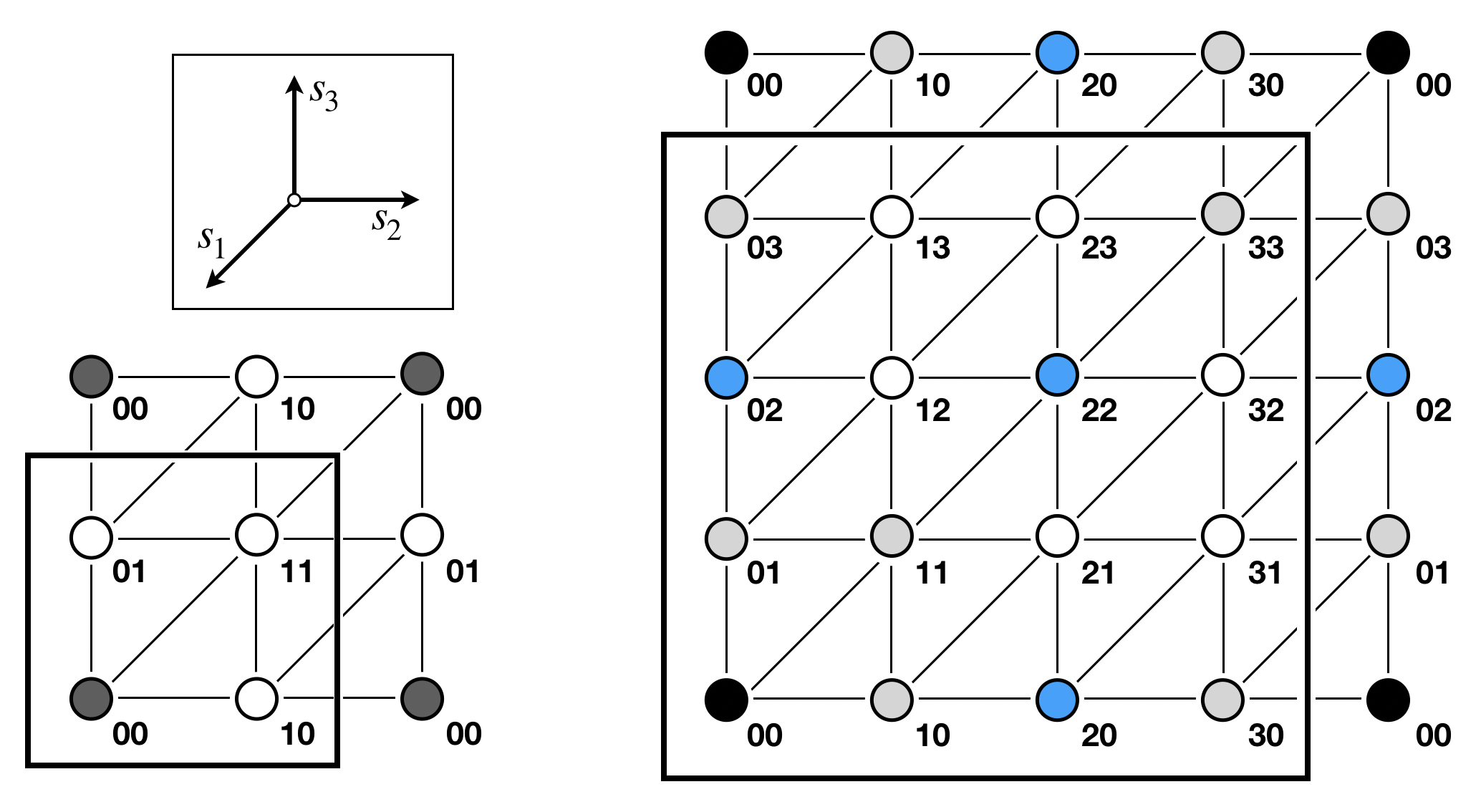}
\caption{Orientation in $\ATo_1$ and $\ATo_2$ -- orientation highlighted in the inset.}
\label{fig:Oriented ATo_n}
\end{figure}
%
%
%
\subsection{Diameter and Antipodals in $\ATo_n$}
\label{subsection:Diameter and Antipodals in ATo_n}
\begin{figure}
\centering
\includegraphics[width=12cm]{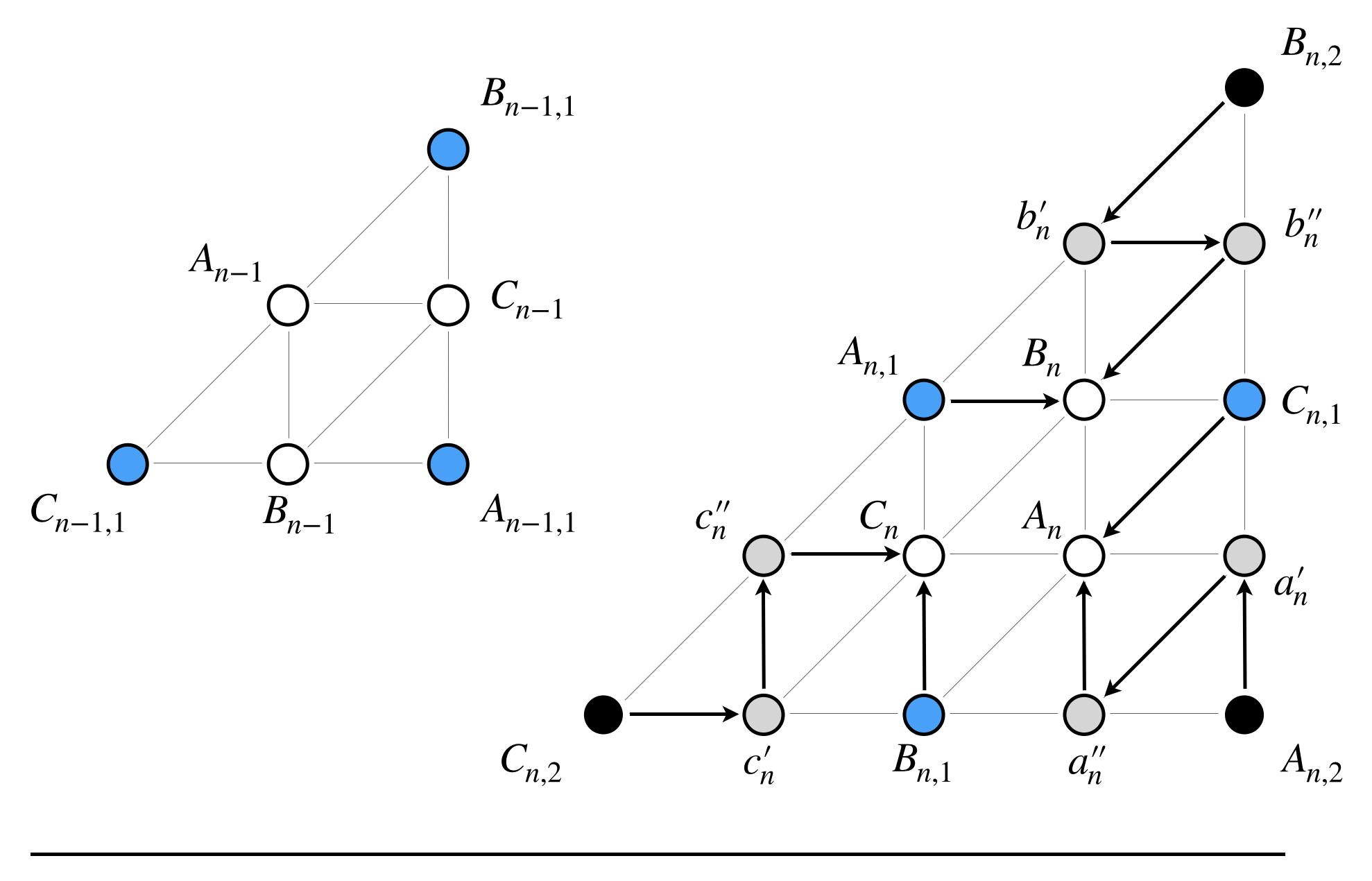}
\includegraphics[width=12cm]{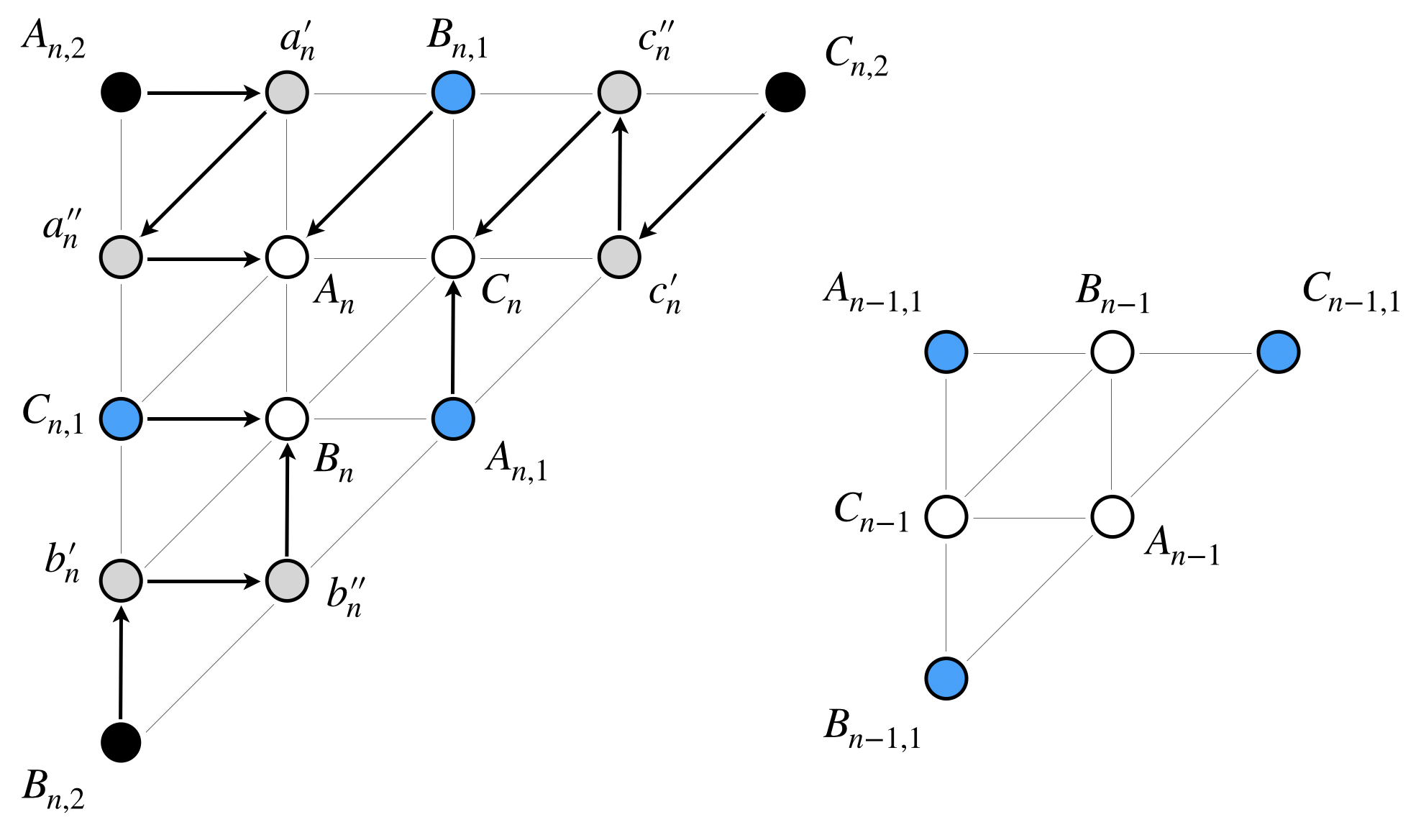}
\caption{Triangular diagrams  in $\ATo_n$ -- for $n$ odd $(\uparrow)$ -- for $n$ even $(\downarrow)$.}
\label{fig:Triangular odd-even diagrams for ATo-n}
\end{figure}

%
%
\noindent We now refer to  digraph $\ATo_n$ in
%
%
Def.\,\ref{definition:  arrowhead and diamond}
%
%
with its generating set $S^+$ highlighted in the inset of 
Fig.\,\ref{fig:Oriented ATo_n}.
Again it appears  as trivial that vertices in subset 
$
	 \{ ( 1,1),( 1,0),( 0,1) \}  
$ 
in  $\ATo_1$ are antipodal and at distance 1 from the origin $ (0,0) $.
We claim that in $\ATo_2$ the subset 
$
	 \{ ( 1,2),( 1,3),( 2,3) \}  
$ 
as well as its symmetric part
$
	 \{ ( 2,1),( 3,1),( 3,2) \}  
$ 
are antipodal at distance 3. To sketch the proof, let us again fix 
$
	x \le y
$
and then extend to the triangular diagrams in 
Fig.\,\ref{fig:Triangular odd-even diagrams for ATo-n}.
%
%
%
%
%
%
\begin{lemma}
\label{lemma:  Oriented diameter in ATo_n}
Let  $\vec{D}_n$ be the oriented diameter of $\ATo_n$. Then 
$
	\vec{D}_0 = 0
$
and
$
	\vec{D}_n =  2 \vec{D}_{n-1} + 1
$
for $ n > 0 $.
\end{lemma}
%
%
\noindent \textit{Proof}:
$\vec{D}_0 = 0$ and in $\ATo_1$ the triple 
$
	( A_1 , B_1 , C_1 ) =  ( ( 1,1),( 1,0),( 0,1) )
$ 
is clearly antipodal and $\vec{D}_1 = 1$.
Let $A_{n,2} $ be the image of $A_{n-1,1}$
whatever the parity of $n$, and assumed from induction at distance
$
	2  ( \vec{D}_{n-1} -1 ).
$
There exists a directed path
\[
	(  \, A_{n,2}  \rightarrow a'_n  \rightarrow a''_n \rightarrow A_n  \, )
\]
of length 3 from $A_{n,2} $ to $A_n$ therefore $A_n$ is at distance
$
	2  ( \vec{D}_{n-1} -1 ) + 3 = 2 \vec{D}_{n-1} + 1 .
$
Idem for 
$
	(  \, B_{n,2}  \rightarrow  B_n  \, )
$
and
$
	(  \, C_{n,2}  \rightarrow  C_n  \, ) .
$

Besides the triple
$
	(\, A_{n,1}, B_{n,1}, C_{n,1} \,)
$
as image of
$
	(\, A_{n-1}, B_{n-1}, C_{n-1} \,)
$
is at distance 
$
	2 \vec{D}_{n-1}
$
by induction. The triple 
$
	( A_n , B_n , C_n )  
$ 
is then also reachable either from
$
	(\, C_{n,1}, A_{n,1}, B_{n,1} \,)
$
for $n$ odd or from
$
	(\, B_{n,1}, C_{n,1}, A_{n,1} \,)
$
for $n$ even. It is therefore antipodal and at distance the diameter
$
	\vec{D}_n =  2 \vec{D}_{n-1} + 1 .
$
$\hfill \Box$
%
%
%
%
\begin{proposition}
\label{proposition: Diameter and antipodals in in ATo_n}
$\ATo_n$ has the diameter
$
	\vec{D}_n =  \sqrt{N}-1 
$
and the number of antipodals
%
%
\begin {itemize}
\item
%
$
    {\mathcal{N_{\vec{A}}}}_{,1} =   \,  \,  \mid \Omega_{1} \cup \overline{\Omega}_{1} \mid \ = \ \mid \Omega_{1} \mid \ = 3
$
\item
$
    {\mathcal{N_{\vec{A}}}}_{,n}  =  \,   \,  \mid \Omega_{n} \cup \overline{\Omega}_{n} \mid \  = 6  \hspace{2cm} ( n > 1 )
$
$\hfill \Box$
\end {itemize}
%
%
\end{proposition}
%
%
%
\subsection{Diameter and Antipodals in $\DTo_n$}
\label{subsection:Diameter and Antipodals in DTo_n}
%
%
\noindent We now refer to  digraph $\DTo_n$ with its generating set 
$
	T^+ = (t_1, t_2, t_3 ) = ( (1,1), (1,0), (0,1) ) 
$
%
in Def.\,\ref{definition:  arrowhead and diamond}.
%
%
\begin{proposition}
\label{proposition: Diameter in DTo_n}
$\DTo_n$ has the diameter
$
	\vec{D}_n =  \sqrt{N}-1 
$
and the number of antipodals
%
%
\begin {itemize}
\item
%
$
    {\mathcal{N_{\vec{D}}}}_{,n}  =  \,   \,   2  \sqrt{N}-1 
$
for any 
$
	n \in  \mathbb{N} .
$
\end {itemize}
%
%
\end{proposition}
%
%
\noindent \textit{Proof}:
The proof here is rather trivial: there is 1 vertex at distance 0, there are 3 vertices at distance 1, there are  $2p + 1$ vertices at distance $p$ and therefore
$
	2 ( 2^n -1 ) + 1 =  2^{n+1} - 1
$
vertices at distance 
$
	2^n - 1 .
$ 
$\hfill \Box$
%
%
%
\section*{Acknowledgement}
%
%
\noindent I am very grateful towards a reviewer of an ancient version of this article and who suggested a complete redefinition of the graphs presented a first time in 
\cite {Deserable-1999},
a redefinition which greatly supported the elaboration of this article.
%
%
%
%

\end{document}